\documentclass[twoside,leqno,10pt, A4]{amsart}
\usepackage{amsfonts}
\usepackage{amsmath}
\usepackage{amscd}
\usepackage{amssymb}
\usepackage{amsthm}
\usepackage{amsrefs}
\usepackage{latexsym}
\usepackage{mathrsfs}
\usepackage{bbm}
\usepackage{amscd}
\usepackage{amssymb}
\usepackage{amsthm}
\usepackage{amsrefs}
\usepackage{latexsym}
\usepackage{mathrsfs}
\usepackage{bbm}
\usepackage{enumerate}
\usepackage{graphicx}
\usepackage{color}
\begin{document}

\newtheorem{theorem}[subsection]{Theorem}
\newtheorem{proposition}[subsection]{Proposition}
\newtheorem{lemma}[subsection]{Lemma}
\newtheorem{corollary}[subsection]{Corollary}
\newtheorem{conjecture}[subsection]{Conjecture}
\newtheorem{prop}[subsection]{Proposition}
\newtheorem{defin}[subsection]{Definition}
\def\phis{\varphi^*}
\newcommand{\Wf}{\mathcal{W}}

\numberwithin{equation}{section}
\newcommand{\mr}{\ensuremath{\mathbb R}}
\newcommand{\mc}{\ensuremath{\mathbb C}}
\newcommand{\dif}{\mathrm{d}}
\newcommand{\intz}{\mathbb{Z}}
\newcommand{\ratq}{\mathbb{Q}}
\newcommand{\natn}{\mathbb{N}}
\newcommand{\comc}{\mathbb{C}}
\newcommand{\rear}{\mathbb{R}}
\newcommand{\prip}{\mathbb{P}}
\newcommand{\uph}{\mathbb{H}}
\newcommand{\fief}{\mathbb{F}}
\newcommand{\majorarc}{\mathfrak{M}}
\newcommand{\minorarc}{\mathfrak{m}}
\newcommand{\sings}{\mathfrak{S}}
\newcommand{\fA}{\ensuremath{\mathfrak A}}
\newcommand{\mn}{\ensuremath{\mathbb N}}
\newcommand{\mq}{\ensuremath{\mathbb Q}}
\newcommand{\half}{\tfrac{1}{2}}
\newcommand{\f}{f\times \chi}
\newcommand{\summ}{\mathop{{\sum}^{\star}}}
\newcommand{\chiq}{\chi \bmod q}
\newcommand{\chidb}{\chi \bmod db}
\newcommand{\chid}{\chi \bmod d}
\newcommand{\sym}{\text{sym}^2}
\newcommand{\hhalf}{\tfrac{1}{2}}
\newcommand{\sumstar}{\sideset{}{^*}\sum}
\newcommand{\sumprime}{\sideset{}{'}\sum}
\newcommand{\sumprimeprime}{\sideset{}{''}\sum}
\newcommand{\sumflat}{\sideset{}{^\flat}\sum}
\newcommand{\shortmod}{\ensuremath{\negthickspace \negthickspace \negthickspace \pmod}}
\newcommand{\V}{V\left(\frac{nm}{q^2}\right)}
\newcommand{\sumi}{\mathop{{\sum}^{\dagger}}}
\newcommand{\mz}{\ensuremath{\mathbb Z}}
\newcommand{\leg}[2]{\left(\frac{#1}{#2}\right)}
\newcommand{\muK}{\mu_{\omega}}
\newcommand{\thalf}{\tfrac12}
\newcommand{\lp}{\left(}
\newcommand{\rp}{\right)}
\newcommand{\Lam}{\Lambda_{[i]}}
\newcommand{\lam}{\lambda}
\newcommand{\af}{\mathfrak{a}}
\newcommand{\sw}{S_{[i]}(X,Y;\Phi,\Psi)}
\newcommand{\lz}{\left(}
\newcommand{\pz}{\right)}
\newcommand{\bfrac}[2]{\lz\frac{#1}{#2}\pz}
\newcommand{\odd}{\mathrm{\ primary}}
\newcommand{\even}{\text{ even}}
\newcommand{\res}{\mathrm{Res}}
\newcommand{\sumn}{\sumstar_{(c,1+i)=1}  w\left( \frac {N(c)}X \right)}
\newcommand{\lab}{\left|}
\newcommand{\rab}{\right|}
\newcommand{\Go}{\Gamma_{o}}
\newcommand{\Ge}{\Gamma_{e}}
\newcommand{\M}{\widehat}
\def\su#1{\sum_{\substack{#1}}}

\theoremstyle{plain}
\newtheorem{conj}{Conjecture}
\newtheorem{remark}[subsection]{Remark}

\newcommand{\pfrac}[2]{\left(\frac{#1}{#2}\right)}
\newcommand{\pmfrac}[2]{\left(\mfrac{#1}{#2}\right)}
\newcommand{\ptfrac}[2]{\left(\tfrac{#1}{#2}\right)}
\newcommand{\pMatrix}[4]{\left(\begin{matrix}#1 & #2 \\ #3 & #4\end{matrix}\right)}
\newcommand{\ppMatrix}[4]{\left(\!\pMatrix{#1}{#2}{#3}{#4}\!\right)}
\renewcommand{\pmatrix}[4]{\left(\begin{smallmatrix}#1 & #2 \\ #3 & #4\end{smallmatrix}\right)}
\def\en{{\mathbf{\,e}}_n}

\newcommand{\ppmod}[1]{\hspace{-0.15cm}\pmod{#1}}
\newcommand{\ccom}[1]{{\color{red}{Chantal: #1}} }
\newcommand{\acom}[1]{{\color{blue}{Alia: #1}} }
\newcommand{\alexcom}[1]{{\color{green}{Alex: #1}} }
\newcommand{\hcom}[1]{{\color{brown}{Hua: #1}} }

\makeatletter
\def\widebreve{\mathpalette\wide@breve}
\def\wide@breve#1#2{\sbox\z@{$#1#2$}%
     \mathop{\vbox{\m@th\ialign{##\crcr
\kern0.08em\brevefill#1{0.8\wd\z@}\crcr\noalign{\nointerlineskip}%
                    $\hss#1#2\hss$\crcr}}}\limits}
\def\brevefill#1#2{$\m@th\sbox\tw@{$#1($}%
  \hss\resizebox{#2}{\wd\tw@}{\rotatebox[origin=c]{90}{\upshape(}}\hss$}
\makeatletter

\title[Bounds for Moments of Dirichlet $L$-functions of fixed modulus on the critical line]{Bounds for Moments of Dirichlet $L$-functions of fixed modulus on the critical line}

\author[P. Gao]{Peng Gao}
\address{School of Mathematical Sciences, Beihang University, Beijing 100191, China}
\email{penggao@buaa.edu.cn}

\author[L. Zhao]{Liangyi Zhao}
\address{School of Mathematics and Statistics, University of New South Wales, Sydney NSW 2052, Australia}
\email{l.zhao@unsw.edu.au}

\begin{abstract}
We study the $2k$-th moment of the family of Dirichlet $L$-functions to a fixed prime modulus on the critical line and establish sharp lower
 bounds for all real $k \geq 0$ and sharp upper bounds for $k$ in the range $0 \leq k \leq 1$.
\end{abstract}

\maketitle



\section{Introduction}\label{sec 1}

  Moments of central values of Dirichlet $L$-functions have been extensively studied in the literature as these values encode significant arithmetic information.  In this paper, we consider the family of Dirichlet $L$-functions to a fixed modulus $q$. To ensure the existence of primitive Dirichlet
  characters modulo $q$, we shall henceforth assume that $q \not \equiv 2 \pmod 4$.  Let $\phis(q)$ denote the number of primitive characters modulo $q$ and $\sum^*$ denote the sum over primitive Dirichlet characters modulo $q$ throughout the paper.
  It is widely believed that (see \cite{R&Sound}) for all real $k  \geq 0$ and explicit constants $C_k$,
\begin{align}
\label{moments}
 \sumstar_{\substack{ \chi \shortmod q }}|L(\tfrac{1}{2},\chi)|^{2k} \sim C_k \phis(q)(\log q)^{k^2}.
\end{align}
 The above was conjectured by J. B. Conrey, D. W. Farmer, J. P. Keating, M. O. Rubinstein and N. C. Snaith in \cite{CFKRS} for all $k \in \natn$.  A. Selberg \cite{Selberg46} had obtained an asymptotic formula concerning a more general twisted second moment of the family of Dirichlet $L$-functions to a fixed modulus, which implies the case of $k=1$ in \eqref{moments}.  Asymptotic evaluation of the case $k=2$ and related results can be found in \cites{HB81, Sound2007, Young2011, BFKMM1, BFKMM, Wu2020, BPRZ}. See also \cites{Sound01,HB2010, Harper, R&Sound1, Radziwill&Sound1, C&L, Gao2024, AC25} for upper and lower bounds of the conjectured order of magnitude.  Using the upper bounds principle developed by  M. Radziwi{\l\l} and K. Soundararajan  \cite{Radziwill&Sound} and the lower bounds principle developed by   W. Heap and K. Soundararajan  \cite{H&Sound}, the first-named author \cite{Gao2024} proved that for large prime $q$ and any real number $k \geq 0$, 
\begin{align*}
\begin{split}
   \sumstar_{\substack{ \chi \shortmod q }}|L(\tfrac{1}{2},\chi)|^{2k} \gg_k & \phis(q)(\log q)^{k^2}, \quad \text{for all } k \geq 0; \\
    \sumstar_{\substack{ \chi \shortmod q }}|L(\tfrac{1}{2},\chi)|^{2k} \ll_k & \phis(q)(\log q)^{k^2}, \quad \text{for } 0 \leq k \leq 1.
\end{split}
\end{align*}  
  
In addition to the average of their central values, there are now growing interests in the evaluation of moments of families of $L$-functions on the critical line. For example, M. Munsch \cite{Munsch17} showed that upper bounds for the shifted moments of the family of Dirichlet $L$-functions to a fixed modulus can be applied to obtain bounds for moments of character sums. The result of Munsch was strengthened by B. Szab\'o \cite{Szab} who proved that under the generalized Riemann hypothesis (GRH), for a large fixed modulus $q$, any positive integer $k$, real tuples ${\bf a} =  (a_1, \ldots, a_k), {\bf t} =  (t_1, \ldots, t_k)$  such that $a_j \geq 0$ and $|t_j| \leq q^A$ for a fixed positive real number $A$,
\begin{align}
\begin{split}
\label{eqn:shiftedMoments}
\sumstar_{\substack{ \chi \shortmod q }}\big| L\big( \tfrac12+it_1,\chi \big) \big|^{2a_1} \cdots \big| L\big( \tfrac12+it_{k},\chi \big) \big|^{2a_{k}} \ll_{\bf{t}, \bf{a}} &  \varphi(q)(\log q)^{a_1^2+\cdots +a_{k}^2} \prod_{1\leq j<l\leq k}  \big|\zeta(1+i(t_j-t_l)+\tfrac 1{\log q}) \big|^{2a_ja_l}, 
\end{split}
\end{align}
 where $\varphi$ denotes the Euler totient function and $\zeta(s)$ the Riemann zeta function. A partial companion lower bounds result can be found in \cite{G&Zhao24-11}, where the authors show, under GRH, that one may replace $\ll_{\bf{t}}$ by  $\gg_{\bf{t}}$ in \eqref{eqn:shiftedMoments} for the case $k=2$ and $q$ being a prime. \newline

  We note that by setting $t_i=t$ in \eqref{eqn:shiftedMoments} and by the estimation that $\zeta(1+\tfrac 1{\log q}) \ll \log q$ (see \cite[Corollary 1.17]{MVa1}), the result of  Szab\'o implies that under GRH,  for any real $k \geq 0$ and $|t| \leq q^A$ for a fixed positive real number $A$,
\begin{align*}
\begin{split}
\sumstar_{\substack{ \chi \shortmod q }}\big| L\big( \tfrac12+it,\chi \big) \big|^{2k} \ll_{t, k} &  \varphi(q)(\log q)^{k^2}. 
\end{split}
\end{align*}  
  
  It is the aim of this paper to establish the above result as well as its companion lower bounds for certain ranges of $k$ and prime $q$ unconditionally.   We begin with our lower bounds. 
\begin{theorem}
\label{thmlowerbound}
   With the notation as above, let $\varepsilon_0$ be fixed with $0<\varepsilon_0<1/4$. Suppose that $q$ is a large prime and $|t| \leq q^{1/4-\varepsilon_0}$. We have for any real number $k \geq 0$,
\begin{align}
\label{lowerbound}
   \sumstar_{\substack{ \chi \shortmod q }}|L(\tfrac{1}{2}+it,\chi)|^{2k} \gg_{t,k} \phis(q)(\log q)^{k^2}.
\end{align}
\end{theorem}

Our upper bounds result is as follows.
\begin{theorem}
\label{thmupperbound}
Under the same assumptions and conditions as Theorem~\ref{thmlowerbound}, we have, for any real number $k$ such that $0 \leq k \leq 1$,
\begin{align}
\label{upperbound}
   \sumstar_{\substack{ \chi \shortmod q }}|L(\tfrac{1}{2}+it,\chi)|^{2k} \ll_{t,k} \phis(q)(\log q)^{k^2}.
\end{align}
\end{theorem}

The combination of Theorems \ref{thmlowerbound} and \ref{thmupperbound} readily gives  the following result concerning the order of magnitude of our family of $L$-functions on the critical line.
\begin{theorem}
\label{thmorderofmag}
Under the same assumptions and conditions as Theorem~\ref{thmlowerbound}, we have, for any real number $k$ such that $0 \leq k \leq 1$,
\begin{align*}
   \sumstar_{\substack{ \chi \shortmod q }}|L(\tfrac{1}{2}+it,\chi)|^{2k} \asymp \phis(q)(\log q)^{k^2}.
\end{align*}
\end{theorem}

   Our proofs of Theorems \ref{thmlowerbound} and \ref{thmupperbound} are obtained by modifying the approaches in \cite{Gao2024}. A key ingredient is the above mentioned result of Selberg on twisted second moment of the family of Dirichlet $L$-functions to a fixed modulus. We remark that in \cite[Theorem 1.3]{Curran24-12}, M. J. Curran obtained result analogous (but more general) to Theorems \ref{thmlowerbound} and \ref{thmupperbound} for the moments of the Riemann zeta function on the critical line.

\section{Proofs of Theorems \ref{thmlowerbound} and \ref{thmupperbound}}
\label{sec 2}

\subsection{Setup}
\label{sec 2'}
 Note first that $\phis(q)=q-2$ when $q$ is prime. As the case $k=0$ is trivial and the case $k=1$ in both \eqref{lowerbound} and \eqref{upperbound} follows from the above-mentioned result of Selberg (see \eqref{Twistedsecmomenttequal} below), we henceforth assume that $k \neq 0,1$.  For two large natural numbers $N$, $M$ depending on $k$ only, we define a sequence of even natural
  numbers $\{ \ell_j \}_{1 \leq j \leq R}$ by $\ell_1= 2\lceil N \log \log q\rceil$ and $\ell_{j+1} = 2 \lceil N \log \ell_j \rceil$ for $j \geq 1$.  Here $R$ denotes the largest natural number satisfying $\ell_R >10^M$.  \newline

   Let ${ P}_1$ denote the set of odd primes not exceeding $q^{1/\ell_1^2}$ and
${ P_j}$ the set of primes lying in the interval $(q^{1/\ell_{j-1}^2}, q^{1/\ell_j^2}]$ for $2\le j\le R$. For $t \in \mr$ and each $1 \leq j \leq R$, we set
\begin{equation*}
{\mathcal P}_j(t, \chi) = \sum_{p\in P_j} \frac{\chi(p)}{p^{1/2+it}}  \quad \mbox{and} \quad {\mathcal Q}_j(t, \chi, k) =\Big (\frac{12 \max(1, k^2)  {\mathcal
P}_j(t, \chi) }{\ell_j}\Big)^{r_k\ell_j},
\end{equation*}
  where $r_k=\lceil 1+1/k \rceil+1$ for $0<k<1$ and $r_k=\lceil k /(2k-1) \rceil+1$ for $k>1$.  We further set ${\mathcal Q}_{R+1}(t, \chi, k)=1$. \newline

  For any non-negative integer $\ell$ and any real number $x$, let
\begin{equation*}
E_{\ell}(x) = \sum_{j=0}^{\ell} \frac{x^{j}}{j!}.
\end{equation*}
  Then we define for each $1 \leq j \leq R$ and real numbers $t, \alpha$,
\begin{align*}
{\mathcal N}_j(t, \chi, \alpha) = E_{\ell_j} (\alpha {\mathcal P}_j(t,\chi)) \quad \mbox{and} \quad \mathcal{N}(t, \chi, \alpha) = \prod_{j=1}^{R} {\mathcal
N}_j(t, \chi,\alpha).
\end{align*}
We follow the convention throughout the paper that the empty product equals to $1$. \newline

  We now proceed as in the proofs of \cite[Lemmas 3.1 and 3.2]{Gao2024} and apply the lower bounds principle of W. Heap and K. Soundararajanand in \cite{H&Sound} and the upper bounds principle of M. Radziwi{\l\l} and K. Soundararajan in \cite{Radziwill&Sound} to obtain the following analogues of those last-mentioned lemmas.
\begin{lemma}
\label{lem1}
 With the notation as above, for $0<k<1$, we have
\begin{align*}
\begin{split}
\sumstar_{\substack{ \chi \shortmod q }} L(\tfrac{1}{2}+it,\chi) & \mathcal{N}(t, \chi, k-1) \mathcal{N}(-t, \overline{\chi}, k) \\
 \leq & \Big ( \sumstar_{\substack{ \chi \shortmod q }}|L(\tfrac{1}{2}+it,\chi)|^{2k} \Big )^{1/2}\Big ( \sumstar_{\substack{ \chi \shortmod q
 }}|L(\tfrac{1}{2}+it,\chi)|^2 |\mathcal{N}(t, \chi, k-1)|^2  \Big)^{(1-k)/2} \\
 & \hspace*{2cm} \times \Big ( \sumstar_{\substack{ \chi \shortmod q }}   \prod^R_{j=1}\big ( |{\mathcal N}_j(t, \chi, k)|^2+ |{\mathcal Q}_j(t, \chi,k)|^2 \big )
 \Big)^{k/2}.
\end{split}
\end{align*}
 For $k>1$, we have
\begin{align}
\label{basicboundkbig}
\begin{split}
\sumstar_{\substack{ \chi \shortmod q }}L(\tfrac{1}{2}+it,\chi) & \mathcal{N}(t, \chi, k-1) \mathcal{N}(-t,\overline{\chi}, k) \\
 \leq & \Big ( \sumstar_{\substack{ \chi \shortmod q }}|L(\tfrac{1}{2}+it,\chi)|^{2k} \Big )^{1/(2k)}\Big ( \sumstar_{\substack{ \chi
 \shortmod q }} \prod^R_{j=1} \big ( |{\mathcal N}_j(t, \chi, k)|^2+ |{\mathcal Q}_j(t,\chi,k)|^2 \big ) \Big)^{(2k-1)/(2k)}.
\end{split}
\end{align}
\end{lemma}

\begin{lemma}
\label{lem2}
 With the notation as above, for $0<k<1$,
\begin{align*}
\begin{split}
\sumstar_{\substack{ \chi \shortmod q }} & |L(\tfrac{1}{2}+it,\chi)|^{2k} \\
&  \ll  \Big ( \sumstar_{\substack{ \chi \shortmod q }}|L(\tfrac{1}{2}+it,\chi)|^2 \sum^{R}_{v=0} \prod^v_{j=1}\Big ( |\mathcal{N}_j(t,\chi, k-1)|^2 \Big ) |{\mathcal
 Q}_{v+1}(t,\chi, k)|^{2}
 \Big)^{k} \\
 & \hspace*{2cm} \times  \Big ( \sumstar_{\substack{ \chi \shortmod q }}  \sum^{R}_{v=0}  \Big (\prod^v_{j=1}|\mathcal{N}_j(t,\chi, k)|^2\Big )|{\mathcal
 Q}_{v+1}(t, \chi, k)|^{2} \Big)^{1-k},
\end{split}
\end{align*}
  where the implied constants depend on $k$ only.
\end{lemma}

The case $k \geq 1$ of Theorem \ref{thmlowerbound} can be obtained using \eqref{basicboundkbig} and applying the similar arguments to those that lead to the proof of the case $0<k < 1$ of Theorem \ref{thmlowerbound}.  Thus we may further assume that $0<k<1$ in the sequel. It then follows from Lemmas \ref{lem1} and \ref{lem2} that in order to establish Theorems \ref{thmlowerbound} and \ref{thmupperbound} for $0<k<1$, it suffices to prove the following three propositions.
\begin{proposition}
\label{Prop4} With the notation as above, let $\varepsilon_0$ be fixed such that $0<\varepsilon_0<1$.  Suppose that $q$ is a large prime and $|t| \leq q^{1/4-\varepsilon_0}$. We have 
\begin{align*}
\sumstar_{\substack{ \chi \shortmod q }}L(\tfrac{1}{2}+it,\chi) \mathcal{N}(t, \chi, k-1)\mathcal{N}(-t, \overline{\chi}, k)  \gg \phis(q)(\log q)^{ k^2
} .
\end{align*}
\end{proposition}

\begin{proposition} \label{Prop5}
Under the same conditions of Proposition~\ref{Prop4}, we have 
\begin{align*}
 \max \Big (  \sumstar_{\substack{ \chi \shortmod q }} & |L(\tfrac{1}{2}+it,\chi)\mathcal{N}(t, \chi, k-1)|^2 ,  \\ 
 & \sumstar_{\substack{ \chi \shortmod q }}|L(\tfrac{1}{2}+it,\chi)|^2 \sum^{R}_{v=0}\Big (\prod^v_{j=1}|\mathcal{N}_j(t, \chi, k-1)|^{2}\Big ) |{\mathcal Q}_{v+1}(t,\chi, k)|^2 \Big ) \ll \phis(q)(\log q)^{ k^2 }.   
\end{align*}
\end{proposition}

\begin{proposition} \label{Prop6}
Under the same conditions of Proposition~\ref{Prop4}, we have  
\begin{align*}
\max \Big ( \sumstar_{\substack{ \chi \shortmod q }}\prod^R_{j=1}\big ( |{\mathcal N}_j(t, \chi, k)|^2+ |{\mathcal Q}_j(t, \chi,k)|^2 \big ),  \sumstar_{\substack{ \chi \shortmod q }} \sum^{R}_{v=0} \Big ( \prod^v_{j=1}|\mathcal{N}_j(t, \chi, k)|^{2}\Big )|{\mathcal
 Q}_{v+1}(t, \chi, k)|^2 \Big )    \ll \phis(q)(\log q)^{ k^2 }.
\end{align*}
\end{proposition}

   As the proof of Proposition \ref{Prop6} is similar to that of \cite[Proposition 3.5]{Gao2024}, we omit it here.  It therefore remains to prove Propositions \ref{Prop4} and \ref{Prop5} in what follows.

\subsection{Twisted second moment of Dirichlet $L$-functions}
 
   As the proof of Proposition \ref{Prop5} relies crucially on the knowledge of twisted second moment of Dirichlet $L$-functions, we explore this topic here. We shall consider a general modulus $q$ in this section and we reserve the letter $p$ for a prime number throughout the paper and we use the notation $\zeta_q(s)$ for the Euler product defining $\zeta(s)$ but omitting those primes dividing $q$, i.e.
 \begin{equation*}
\zeta_q(s)=\prod_{p\mid q}\left(1-\frac{1}{p^s}\right)\zeta(s).
\end{equation*}
 In \cite{Selberg46} (see also \cite[Section 4]{Conrey07}), A. Selberg proved that for positive integers $h$ and $k$ with $(h,k)=(hk, q)=1$, one has, for $s=\sigma+it$ and $s'=\sigma'+it'$ and $0<\sigma$, $\sigma' <1$,
\begin{align*}
\begin{split}
\sumstar_{\chi \shortmod q} L(s,\chi)L(s',\overline{\chi}) \chi(h)\overline{\chi}(k) 
 = \frac{\varphi(q)}{h^{s'}k^{s}}\zeta_q(s+s') & +\frac{\varphi(q)^2(2\pi)^{s+s'-1}}{\pi q^{s+s'}h^{1-s}k^{1-s'}}
\Gamma(1-s)\Gamma(1-s')\cos \tfrac \pi 2 (s-s')\zeta(2-s-s')\\
&+O\left(\frac{|ss'|}{\sigma\sigma'(1-\sigma)(1-\sigma')}
\big(q^{\varepsilon}(hq^{1-\sigma}+kq^{1-\sigma'}+hkq^{1-\sigma-\sigma'})\big)\right). 
\end{split}
\end{align*}
  We deduce from the above upon setting $s=1/2+it, s'=1/2+it'$ with $t, t'\in \mr$ that for $(h,k)=(hk, q)=1$,
\begin{align}
\label{Twistedsecmomentgen}
\begin{split}
 \sumstar_{\chi \shortmod q} & L(\tfrac12+it,\chi)L(\tfrac12+it',\overline{\chi}) \chi(h)\overline{\chi}(k) \\
 =& \frac{\varphi(q)}{h^{1/2+it'}k^{1/2+it}}\zeta_q(1+i(t+t'))  +\frac{\varphi(q)^2 (2\pi)^{i(t+t')}}{\pi q^{1+i(t+t')}h^{1/2-it}k^{1/2-it'}}
\Gamma(\tfrac12-it)\Gamma(\tfrac12-it')\cos \tfrac \pi 2i (t-t')\zeta(1-i(t+t'))\\
& \hspace*{3cm} +O\left((|t|+1)(|t'|+1)((h+k)q^{1/2+\varepsilon}+hkq^{\varepsilon})\right). 
\end{split}
\end{align}   
   We now consider the case $t \neq 0$ and $t'\rightarrow -t$ in \eqref{Twistedsecmomentgen}. From \cite[Corollary 1.17]{MVa1},
\begin{align*}
\begin{split}
 \zeta(1+it)=\frac 1{it}+O(1), \quad |t| \leq 1.
\end{split}
\end{align*}   
  Also, it follows from \cite[p. 73, (3)]{Da} that
\begin{align*}
\begin{split}
 \Gamma(\tfrac12-it)\Gamma(\tfrac12+it)=\frac {\pi}{\sin (\pi (\tfrac12+i t))}=\frac {\pi}{\cos (\pi i t)}. 
\end{split}
\end{align*}      
   
   The above formulas gives that the right-hand side of \eqref{Twistedsecmomentgen} remains holomorphic in the process that $t \neq 0$ and $t'\rightarrow t$. Note moreover that
\begin{align*}
\begin{split}   
   (2\pi)^{i(t+t')}=& 1+\log (2\pi) i(t+t')+O((t+t')^2), \\
   q^{-i(t+t')}= & 1-(\log q) i(t+t')+O((t+t')^2), \\ 
   \prod_{p\mid q}\left(1-\frac{1}{p^{1+i(t+t')}}\right) =&   \prod_{p\mid q}\left(1-\frac{1}{p}\right)+\prod_{p\mid q}\left(1-\frac{1}{p}\right)\sum_{p|q}\frac {i(\log p)(t+t')}{p-1}+O((t+t')^2), \\
   h^{-i(t+t')}=& 1-i(\log h) (t+t')+O((\log h)^2(t+t')^2), \\
   k^{i(t+t')}=& 1+i(\log k) (t+t')+O((\log k)^2(t+t')^2), \\ 
   \Gamma(\tfrac12+it-i(t+t'))=&  \Gamma(\tfrac12+it)-i\Gamma'(\tfrac12+it)(t+t')+O((t+t')^2), \\
   \cos \tfrac \pi 2i (t-t')=&  \cos (i\pi t)+\tfrac \pi 2i \sin (i\pi t)(t+t')+O((t+t')^2).
\end{split}
\end{align*}   
  We then conclude that
\begin{align}
\label{Twistedsecmomenttequal}
\begin{split}
 \sumstar_{\chi \shortmod q} & \big |L(\tfrac12+it,\chi)\big |^2 \chi(h)\overline{\chi}(k) 
= \lim_{t'\rightarrow -t}\sumstar_{\chi \bmod q} L(\tfrac12+it,\chi)L(\tfrac12+it',\overline{\chi}) \chi(h)\overline{\chi}(k) \\
 =& \frac{\varphi(q)}{h^{1/2-it}k^{1/2+it}}\prod_{p\mid q}\left(1-\frac{1}{p}\right)\sum_{p|q}\frac {\log p}{p-1}-\frac{\varphi(q)}{ h^{1/2-it}k^{1/2+it}}\log h+\frac{\varphi(q)^2}{qh^{1/2-it}k^{1/2+it}}\log \Big( \frac {q}{2\pi} \Big)\\
 & \hspace*{2cm} +\frac{\varphi(q)^2}{qh^{1/2-it}k^{1/2+it}}\frac {\Gamma'(\tfrac12+it)}{\Gamma(\tfrac12+it)}-\frac{\varphi(q)^2}{qh^{1/2-it}k^{1/2+it}}\tfrac \pi 2\tan(i\pi t)-\frac{\varphi(q)^2}{qh^{1/2-it}k^{1/2+it}}\log k \\
 & \hspace*{2cm}  +O\left((|t|+1)^2((h+k)q^{1/2+\varepsilon}+hkq^{\varepsilon})\right). 
\end{split}
\end{align}

\subsection{Proof of Proposition \ref{Prop4}}
\label{sec 4}

   As shown in the proof of \cite[Proposition 3.3]{Gao2024}, one sees that by taking $M$ large enough, ${\mathcal N}(t, \chi, k-1)$ and ${\mathcal N}(-t, \overline \chi, k)$  are short Dirichlet
    polynomials with lengths at most $q^{2/10^{M}}$. Also, we have
\begin{align}
\label{prodNbound}
 & {\mathcal N}(t, \chi, k-1){\mathcal N}(-t, \overline \chi, k) \ll q^{4/10^{M}}.
\end{align}

Moreover it is shown in \cite{R&Sound} that for $X \geq 1$,
\begin{equation} \label{approxfneq}
L(\half+it, \chi)=\sum_{m \leq X}\frac {\chi(m)}{m^{1/2+it}}+O\left( \frac {(|t|+1)\sqrt{q}\log q}{\sqrt{X}} \right).
\end{equation}

  We now write for simplicity, 
\begin{align*}
 {\mathcal N}(t, \chi, k-1)= \sum_{a  \leq q^{2/10^{M}}} \frac{x_a}{a^{1/2+it}} \chi(a) \quad \mbox{and} \quad \mathcal{N}(-t, \overline{\chi}, k) = \sum_{b  \leq
 q^{2/10^{M}}} \frac{y_b}{b^{1/2-it}}\overline{\chi}(b).
\end{align*}
  
  We then deduce from \eqref{prodNbound} and \eqref{approxfneq} that
\begin{align*}
 \sumstar_{\substack{ \chi \shortmod q }} & L(\tfrac{1}{2}+it,\chi)\mathcal{N}(t, \chi, k-1) \mathcal{N}(-t, \overline{\chi}, k)  \\
=&  \sumstar_{\substack{ \chi \shortmod q }}\sum_{m \leq X}\frac {\chi(m)}{m^{1/2+it}}\mathcal{N}(t, \chi,
k-1)\mathcal{N}(-t, \overline{\chi}, k)+O\left( \frac {(|t|+1)\phis(q)q^{1/2+4/10^{M}}\log q}{\sqrt{X}} \right) \\
=& \sum_{\substack{ \chi \shortmod q }}\sum_{m \leq X}\frac {\chi(m)}{m^{1/2+it}}\mathcal{N}(t, \chi,
k-1)\mathcal{N}(-t, \overline{\chi}, k)+O\left( \sqrt{X}q^{4/10^{M}}+\frac {(|t|+1)\phis(q)q^{1/2+4/10^{M}}\log q}{\sqrt{X}} \right)  \\
=& \phis(q) \sum_{a} \sum_{b} \sum_{\substack{m \leq X \\ am \equiv b \bmod q}}\frac {x_a y_b}{(am)^{1/2+it}b^{1/2-it}}+O\left( \sqrt{X}q^{4/10^{M}}+\frac {(|t|+1)\phis(q)q^{1/2+4/10^{M}}\log q}{\sqrt{X}} \right) .
\end{align*}
  
  Similar to the proof of \cite[Proposition 3.3]{Gao2024}, the contribution from the terms $am=b+l q$ with $l \geq 1$ in the last expression above (note that as $b <q$, we can not have $b > am$ in our case) is $\ll  \sqrt{X}q^{2/10^{M}}$.  \newline

   We now set $X=(|t|+1)q^{3/2+4/10^{M}}$ to see that by taking $M$ large enough, the contributions from various error terms above when $|t| \leq q^{1/4-\varepsilon_0}$ can be ignored. We thus deduce that
\begin{align*}
& \sumstar_{\substack{ \chi \shortmod q }}L(\tfrac{1}{2}+it,\chi)\mathcal{N}(t, \chi, k-1) \mathcal{N}(-t, \overline{\chi}, k) 
\gg \phis(q) \sum_{a} \sum_{b} \sum_{\substack{m \leq X \\ am = b }}\frac {x_a y_b}{\sqrt{abm}}. 
\end{align*}
  Note that the right-hand side expression above is independent of $t$. We then proceed as in the proof of \cite[Proposition 3.3]{Gao2024} to see that the right-hand side expression above is $\gg \phis(q)(\log q)^{k^2}$.
 The assertion of Proposition \ref{Prop4} now follows from this. This completes the proof.

\subsection{Proof of Proposition \ref{Prop5}}
\label{sec 5}

   As the proofs are similar, it suffices to focus on the second term in the maximum and show that
\begin{align}
\label{sumLsquareNQ}
 \sumstar_{\substack{ \chi \shortmod q }}|L(\tfrac{1}{2}+it,\chi)|^2 \sum^{R}_{v=0}\Big (\prod^v_{j=1}|\mathcal{N}_j(t, \chi, k-1)|^{2}\Big ) |{\mathcal
 Q}_{v+1}(t,\chi, k)|^2  
  \ll &
\phis(q)(\log q)^{ k^2 }.   
\end{align}   
   
   Let $p_{v+1}(n)$ be the function whose range is $\{0,1\}$, and that $p_{v+1}(n)=1$ if and only if $n$ is composed of exactly $r_k\ell_{v+1}$ primes (counted with multiplicity), all from the interval $P_{v+1}$. Together with the notations in Section \ref{sec 4}, we now write that
\begin{align*}
  {\mathcal P}_{v+1}(t, \chi)^{r_k\ell_{v+1}} =&  \sum_{ \substack{ n_{v+1}}} \frac{1}{\sqrt{n_{v+1}}}\frac{(r_k\ell_{v+1})!
  }{w(n_{v+1})}\chi(n_{v+1})p_{v+1}(n_{v+1}).
\end{align*}
Recall that $r_k=\lceil 1+1/k \rceil+1$ when $0<k<1$. This allows us to further write 
\begin{align}
\label{ProdNQ}
\begin{split}
 \Big (\prod^v_{j=1}|\mathcal{N}_j(t, \chi, k-1)|^2 \Big )|{\mathcal
 Q}_{v+1}(t, \chi, k)|^{2} = \Big( \frac{12  }{\ell_{v+1}}\Big)^{2r_k\ell_{v+1}}((r_k\ell_{v+1})!)^2 \sum_{a,b \leq q^{2r_k/10^{M}}} \frac{u_a u_b}{a^{1/2+it}b^{1/2-it}}\chi(a)\overline{\chi}(b).
\end{split}
\end{align}
 Here we notice that $\prod^v_{j=1}\mathcal{N}_j(t, \chi, k-1) \cdot {\mathcal
 Q}_{v+1}(t,\chi, k)$ is a short Dirichlet polynomial whose length is at most $q^{2r_k/10^{M}}$. 
Note also that $u_a, u_b \leq 1$ for all $a$, $b$. \newline

 We deduce from \eqref{ProdNQ} that
\begin{align*}
\begin{split}
\sumstar_{\substack{ \chi \shortmod q }} & |L(\tfrac12+it, \chi)|^2\Big (\prod^v_{j=1}|\mathcal{N}_j(t,\chi, k-1)|^2 \Big )|{\mathcal
 Q}_{v+1}(t,\chi, k)|^{2} \\
=& \Big( \frac{12  }{\ell_{v+1}}\Big)^{2r_k\ell_{v+1}}((r_k\ell_{v+1})!)^2 \sum_{a,b \leq q^{2r_k/10^{M}}} \frac{u_a u_b}{a^{1/2+it}b^{1/2-it}} \sumstar_{\substack{ \chi \shortmod q }} |L(\tfrac12+it, \chi)|^2\chi(a) \overline{\chi}(b).
\end{split}
\end{align*}

With $(a,b)$ denoting the greatest common divisor of $a$ and $b$, we write $a=(a,b)\cdot a/(a,b), b=(a,b)\cdot b/(a,b)$, leading to $\chi(a) \overline{\chi}(b)=\chi(a/(a,b)) \overline{\chi}(b/(a,b))$. Note that $(a/(a,b), b/(a,b))=1$. Moreover, as both $a, b \leq q^{2r_k/10^{M}}<q$ and $q$ is a prime, we have $(a/(a,b)\cdot b/(a,b), q)=1$. Thus, we are able to apply \eqref{Twistedsecmomenttequal} to evaluate the inner sum, arriving at
\begin{align}
\label{LNsquaresum}
\begin{split}
\sumstar_{\substack{ \chi \shortmod q }} & |L(1/2+it, \chi)|^2\Big (\prod^v_{j=1}|\mathcal{N}_j(t,\chi, k-1)|^2 \Big )|{\mathcal
 Q}_{v+1}(t,\chi, k)|^{2} \\
=& \Big( \frac{12  }{\ell_{v+1}}\Big)^{2r_k\ell_{v+1}}((r_k\ell_{v+1})!)^2 \sum_{a,b \leq q^{2r_k/10^{M}}} \frac{(a,b)u_a u_b}{ab} (S_1+S_2)\\
&\hspace*{2cm} + O\left( \Big( \frac{12  }{\ell_{v+1}}\Big)^{2r_k\ell_{v+1}}((r_k\ell_{v+1})!)^2(|t|+1)^2\sum_{a,b \leq q^{2r_k/10^{M}}} \frac{1}{\sqrt{ab}} ((a+b)q^{1/2+\varepsilon}+abq^{\varepsilon})  \right),
\end{split}
\end{align}
 where 
\begin{align*}
\begin{split}
 S_1 =& \varphi(q)\prod_{p\mid q}\left(1-\frac{1}{p}\right)\sum_{p|q}\frac {\log p}{p-1}+\frac{\varphi(q)^2}{q}\log (\frac {q}{2\pi})- \varphi(q)\log \frac {a}{(a,b)}- \frac {\varphi(q)^2}{q}\log \frac {b}{(a,b)} \quad \mbox{and} \\
 S_2=&\frac {\varphi(q)^2}{q}\frac {\Gamma'(\tfrac12+it)}{\Gamma(\tfrac12+it)}-\frac{\varphi(q)^2}{q}\tfrac \pi 2\tan(i\pi t). 
\end{split}
\end{align*}   

 Note that by \cite[(5.3)]{Gao2024}, we have
\begin{align*}
 & \Big( \frac{12 }{\ell_{v+1}} \Big)^{2r_k\ell_{v+1}}((r_k \ell_{v+1})!)^2
\ll q^{\varepsilon}.
\end{align*}

   It follows from this that upon taking $M$ large enough, the $O$-term in \eqref{LNsquaresum} for $|t| \leq q^{1/4-\varepsilon_0}$ is
$ \ll q^{1-\varepsilon}$. Note again that the main term involving with $S_1$ in \eqref{LNsquaresum} is independent of $t$.  We then proceed as in the proof of \cite[Proposition 3.4]{Gao2024} to see that 
\begin{align}
\label{S1sum}
\begin{split}
 \Big( \frac{12  }{\ell_{v+1}}\Big)^{2r_k\ell_{v+1}}((r_k\ell_{v+1})!)^2 \sum_{a,b \leq q^{2r_k/10^{M}}} \frac{(a,b)u_a u_b}{ab}S_1 \ll \phis(q) e^{-\ell_{v+1}/2}(\log q)^{k^2}.
\end{split}
\end{align}
  Due to the presence of the term $\frac{\varphi(q)^2}{q}\log (\frac {q}{2\pi})$ in $S_1$, the above further implies that  
\begin{align}
\label{absum}
\begin{split}
 \frac{\varphi(q)^2}{q}\Big( \frac{12  }{\ell_{v+1}}\Big)^{2r_k\ell_{v+1}}((r_k\ell_{v+1})!)^2 \sum_{a,b \leq q^{2r_k/10^{M}}} \frac{(a,b)u_a u_b}{ab} \ll \phis(q) e^{-\ell_{v+1}/2}(\log q)^{k^2-1}.
\end{split}
\end{align}

  Note by \cite[Theorem C.1]{MVa1}, we have for $|t| \leq q^{1/4}$, 
\begin{align*}
\begin{split}
 \frac {\Gamma'(\tfrac12+it)}{\Gamma(\tfrac12+it)} \ll \log (1+|t|) \ll \log q. 
\end{split}
\end{align*}  
  As we have $\tan(i\pi t) \ll 1$, we deduce from the above and \eqref{absum} that 
\begin{align}
\label{S2sum}
\begin{split}
 \Big( \frac{12  }{\ell_{v+1}}\Big)^{2r_k\ell_{v+1}}((r_k\ell_{v+1})!)^2 \sum_{a,b \leq q^{2r_k/10^{M}}} \frac{(a,b)u_a u_b}{ab}S_2 \ll \phis(q) e^{-\ell_{v+1}/2}(\log q)^{k^2}.
\end{split}
\end{align}

 Thus we deduce from \eqref{LNsquaresum}, \eqref{S1sum} and \eqref{S2sum} that
\begin{align*}
\begin{split}
  \sumstar_{\substack{ \chi \shortmod q }}|L(\tfrac{1}{2}+it,\chi)|^2 \Big (\prod^v_{j=1}|\mathcal{N}_j(t, \chi, k-1)|^2 \Big )|{\mathcal
 Q}_{v+1}(t,\chi, k)|^{2} \ll \phis(q)e^{-\ell_{v+1}/2}(\log q)^{k^2}.
\end{split}
\end{align*}
 As the sum over $e^{-\ell_j/2}$ converges, the bound in \eqref{sumLsquareNQ} follows from the above. This completes the proof of the proposition.

\vspace*{.5cm}

\noindent{\bf Acknowledgments.}  P. G. is supported in part by NSFC grant 12471003 and L. Z. by the FRG Grant PS71536 at the University of New South Wales.

\bibliography{biblio}

\begin{bibdiv}
\begin{biblist}

\bib{AC25}{article}{
      author={Arguin, L.-P.},
      author={Creighton, N.},
       title={Upper bounds on large deviations of {D}irichlet {$L$}-functions
  in the $q$-aspect},
        date={to appear},
     journal={J. Number Theory},
}

\bib{BFKMM}{article}{
      author={Blomer, V.},
      author={Fouvry, \'{E}.},
      author={Kowalski, E.},
      author={Michel, P.},
      author={Mili\'{c}evi\'{c}, D.},
       title={On moments of twisted {$L$}-functions},
        date={2017},
        ISSN={0002-9327},
     journal={Amer. J. Math.},
      volume={139},
      number={3},
       pages={707\ndash 768},
}

\bib{BFKMM1}{article}{
      author={Blomer, V.},
      author={Fouvry, \'{E}.},
      author={Kowalski, E.},
      author={Michel, P.},
      author={Mili\'{c}evi\'{c}, D.},
       title={Some applications of smooth bilinear forms with {K}loosterman
  sums},
        date={2017},
     journal={Tr. Mat. Inst. Steklova},
      volume={296},
       pages={Analiticheskaya i Kombinatornaya Teoriya Chisel, 24\ndash 35},
        note={English version published in Proc. Steklov Inst. Math.
  {{\bf{2}}96} (2017), no. 1, 18--29},
}

\bib{BPRZ}{article}{
      author={Bui, H.~M.},
      author={Pratt, K.},
      author={Robles, N.},
      author={Zaharescu, A.},
       title={Breaking the {$\frac12$}-barrier for the twisted second moment of
  {D}irichlet {$L$}-functions},
        date={2020},
     journal={Adv. Math.},
      volume={370},
       pages={107175, 40pp},
}

\bib{C&L}{article}{
      author={Chandee, V.},
      author={Li, X.},
       title={Lower bounds for small fractional moments of {D}irichlet
  {$L$}-functions},
        date={2013},
     journal={Int. Math. Res. Not. IMRN},
      number={19},
       pages={4349\ndash 4381},
}

\bib{Conrey07}{article}{
      author={Conrey, J.~B.},
       title={The mean-square of {D}irichlet {$L$}-functions},
        date={Preprint},
        note={arXiv:0708.2699},
}

\bib{CFKRS}{article}{
      author={Conrey, J.~B.},
      author={Farmer, D.~W.},
      author={Keating, J.~P.},
      author={Rubinstein, M.~O.},
      author={Snaith, N.~C.},
       title={Integral moments of {$L$}-functions},
        date={2005},
     journal={Proc. London Math. Soc. (3)},
      volume={91},
      number={1},
       pages={33\ndash 104},
}

\bib{Curran24-12}{article}{
      author={Curran, M.~J.},
       title={Freezing transition and moments of moments of the {R}iemann zeta
  function},
        date={2024},
     journal={Q. J. Math.},
      volume={75},
      number={4},
       pages={1481–1505},
}

\bib{Da}{book}{
      author={Davenport, H.},
       title={{M}ultiplicative {N}umber {T}heory},
     edition={Third edition},
      series={Graduate Texts in Mathematics},
   publisher={Springer-Verlag},
     address={Berlin, etc.},
        date={2000},
      volume={74},
}

\bib{Gao2024}{article}{
      author={Gao, P.},
       title={Bounds for moments of {D}irichlet {$L$}-functions to a fixed
  modulus},
        date={2024},
     journal={Math. Z.},
      volume={307},
      number={4},
       pages={Paper No. 70, 18pp},
}

\bib{G&Zhao24-11}{article}{
      author={Gao, P.},
      author={Zhao, L.},
       title={Lower bounds for shifted moments of {D}irichlet {$L$}-functions
  of fixed modulus},
        date={Preprint},
        note={arXiv:2411.03692},
}

\bib{Harper}{article}{
      author={Harper, A.~J.},
       title={Sharp conditional bounds for moments of the {R}iemann zeta
  function},
        date={Preprint},
        note={arXiv:1305.4618},
}

\bib{H&Sound}{article}{
      author={Heap, W.},
      author={Soundararajan, K.},
       title={Lower bounds for moments of zeta and {$L$}-functions revisited},
        date={2022},
     journal={Mathematika},
      volume={68},
      number={1},
       pages={1\ndash 14},
}

\bib{HB81}{article}{
      author={Heath-Brown, D.~R.},
       title={The fourth power mean of {D}irichlet's {$L$}-functions},
        date={1981},
     journal={Analysis},
      volume={1},
      number={1},
       pages={25\ndash 32},
}

\bib{HB2010}{article}{
      author={Heath-Brown, D.~R.},
       title={Fractional moments of {D}irichlet {$L$}-functions},
        date={2010},
     journal={Acta Arith.},
      volume={145},
      number={4},
       pages={397\ndash 409},
}

\bib{MVa1}{book}{
      author={Montgomery, H.~L.},
      author={Vaughan, R.~C.},
       title={{M}ultiplicative number theory. {I}. {C}lassical theory},
      series={Cambridge Studies in Advanced Mathematics},
   publisher={Cambridge University Press},
     address={Cambridge},
        date={2007},
      volume={97},
}

\bib{Munsch17}{article}{
      author={Munsch, M.},
       title={Shifted moments of {$L$}-functions and moments of theta
  functions},
        date={2017},
     journal={Mathematika},
      volume={63},
      number={1},
       pages={196\ndash 212},
}

\bib{Radziwill&Sound1}{article}{
      author={Radziwi{\l\l}, M.},
      author={Soundararajan, K.},
       title={Continuous lower bounds for moments of zeta and {$L$}-functions},
        date={2013},
     journal={Mathematika},
      volume={59},
      number={1},
       pages={119\ndash 128},
}

\bib{Radziwill&Sound}{article}{
      author={Radziwi{\l \l}, M.},
      author={Soundararajan, K.},
       title={Moments and distribution of central {$L$}-values of quadratic
  twists of elliptic curves},
        date={2015},
     journal={Invent. Math.},
      volume={202},
      number={3},
       pages={1029\ndash 1068},
}

\bib{R&Sound1}{article}{
      author={Rudnick, Z.},
      author={Soundararajan, K.},
       title={Lower bounds for moments of {$L$}-functions: symplectic and
  orthogonal examples},
     journal={in: Multiple {D}irichlet series, automorphic forms, and analytic
  number theory, 293--303, Proc. Sympos. Pure Math.},
      volume={75},
       pages={Amer. Math. Soc., Providence, RI, 2006.},
}

\bib{R&Sound}{article}{
      author={Rudnick, Z.},
      author={Soundararajan, K.},
       title={Lower bounds for moments of {$L$}-functions},
        date={2005},
     journal={Proc. Natl. Acad. Sci. USA},
      volume={102},
      number={19},
       pages={6837\ndash 6838},
}

\bib{Selberg46}{article}{
      author={Selberg, A.},
       title={Contributions to the theory of {D}irichlet's {$L$}-functions},
        date={1946},
     journal={Skr. Norske Vid.-Akad. Oslo I},
      volume={1946},
      number={3},
       pages={62},
}

\bib{Sound2007}{article}{
      author={Soundararajan, K.},
       title={The fourth moment of {D}irichlet {$L$}-functions},
     journal={in: Analytic number theory, 239–246, Clay Math. Proc.,},
      volume={7},
       pages={Amer. Math. Soc., Providence, RI, 2007.},
}

\bib{Sound01}{article}{
      author={Soundararajan, K.},
       title={Moments of the {R}iemann zeta function},
        date={2009},
     journal={Ann. of Math. (2)},
      volume={170},
      number={2},
       pages={981\ndash 993},
}

\bib{Szab}{article}{
      author={Szab\'o, B.},
       title={High moments of theta functions and character sums},
        date={2024},
     journal={Mathematika},
      volume={70},
      number={2},
       pages={Paper No. e12242, 37 pp.},
}

\bib{Wu2020}{article}{
      author={Wu, X.},
       title={The fourth moment of {D}irichlet {$L$}-functions at the central
  value},
        date={2023},
     journal={Math. Ann.},
      volume={387},
      number={3-4},
       pages={1199\ndash 1248},
}

\bib{Young2011}{article}{
      author={Young, M.~P.},
       title={The fourth moment of {D}irichlet {$L$}-functions},
        date={2011},
     journal={Ann. of Math. (2)},
      volume={173},
      number={1},
       pages={1\ndash 50},
}

\end{biblist}
\end{bibdiv}
\bibliographystyle{amsxport}

\end{document}